\documentclass[proceedings,submission]{amsart}
\usepackage[latin1]{inputenc}
\usepackage[T1]{fontenc}
\usepackage{ae,aecompl}

\usepackage{subfigure}

\usepackage{amssymb}
\usepackage{amsmath}
\usepackage{units}

\usepackage[english]{babel}

\usepackage[normalem]{ulem}

\usepackage{tikz,ifthen}
\usetikzlibrary{matrix,shapes,arrows,chains}
\usetikzlibrary{calc}
\usetikzlibrary{patterns}
\usetikzlibrary{positioning} 

\definecolor{darkgreen}{rgb}{0,0.7,0}
\definecolor{darkred}{rgb}{0.7,0,0}
\definecolor{darkblue}{rgb}{0,0,0.7}
\usepackage{pdflscape}

\usepackage{enumerate}
\usepackage{verbatim}
\usepackage{hyperref}

\usepackage{graphicx}
\usepackage{rotating}
\usepackage{algorithmic}
\usepackage{algorithm}
\usepackage{xspace}

\usepackage{ifpdf}
\ifpdf
\fi

\newtheorem{theorem}{Theorem}[section]

\newtheorem{proposition}[theorem]{Proposition}
\newtheorem{corollary}[theorem]{Corollary}
\newtheorem{definition}[theorem]{Definition}

\newtheorem{problem}[theorem]{Problem}

\newcommand{\x}{\mathbf{x}}
\renewcommand{\a}{\mathbf{a}}
\renewcommand{\v}{\mathbf{v}}

\newcommand{\orbsum}[1][G]{\mathfrak{\sum_{\operatorname{Orb}(#1)}}}
\newcommand{\suchthat}{\mid}

\newcommand{\sg}[1][n]{\mathfrak{S}_{#1}}
\newcommand{\sgs}[1][G]{sgs({#1})}

\newcommand{\NN}{\mathbb{N}}

\newcommand{\KK}{\mathbb{K}}

\newcommand{\Sage}{\texttt{Sage}\xspace}
\newcommand{\sagecombinat}{\texttt{Sage-Combinat}\xspace}
\newcommand{\gap}{\texttt{GAP}\xspace}

\newcommand{\cython}{\texttt{Cython}\xspace}

\newcommand{\ind}{\hspace{4ex}}

\definecolor{sagecolor}{rgb}{.78, .36, .04}
\newcommand{\sageex}[1]{\vspace{0.5em} {\small \tt #1} \vspace{0.5em}}
\def\sageret#1{\begin{center} #1 \end{center}}

\def\sagepromt{{\normalsize \color{sagecolor} sage: }}

\usepackage{ifthen}
\newboolean{draft}
\setboolean{draft}{true}
\ifdraft
\newcommand{\TODO}[2][To do: ]{\textcolor{red}{\textbf{#1#2}}}
\newcommand{\INFO}[2][Info: ]{\textcolor{red}{\textbf{#1#2}}}
\else
\newcommand{\TODO}[2][]{}
\newcommand{\INFO}[2][]{}
\fi

\title[Generation modulo the action of a permutation
  group]{Generating tuples of integers modulo the action of a
  permutation group and applications}

\author{Nicolas Borie}

\address{Univ. Paris Est
  Marne-La-Vall\'ee, Laboratoire d'Informatique Gaspard Monge, Cit\'e
  Descartes, B\^at Copernic -- 5, bd Descartes Champs sur Marne 77454
  Marne-la-Vall\'ee Cedex 2, France}

\begin{document}


\maketitle

\begin{abstract}
  Originally motivated by algebraic invariant theory, we present an
  algorithm to enumerate integer vectors modulo the action of a
  permutation group. This problem generalizes the generation of
  unlabeled graph up to an isomorphism. In this paper, we present the
  full development of a generation engine by describing the
  related theory, establishing a mathematical and practical
  complexity, and exposing some benchmarks. We next show two
  applications to effective invariant theory and effective Galois
  theory.

  Initialement motiv\'e par la th\'eorie alg\'ebrique des invariants,
  nous pr\'esentons une strat\'egie algorithmique pour \'enum\'erer
  les vecteurs d'entiers modulo l'action d'un groupe de
  permutations. Ce probl\`eme g\'en\'eralise le probl\`eme
  d'\'enum\'eration des graphes non \'etiquet\'es. Dans cet article,
  nous d\'eveloppons un moteur complet d'\'enum\'eration en expliquant
  la th\'eorie sous-jacente, nous \'etablissons des bornes de
  complexit\'e pratiques et th\'eoriques et exposons quelques bancs
  d'essais. Nous d\'etaillons ensuite deux applications
  th\'eoriques en th\'eorie effective des invariants et en th\'eorie
  de Galois effective.

\end{abstract}


\section{Introduction}

Let $G$ be a group of permutations, that is, a subgroup of some
symmetric group $\sg[n]$. Several problems in effective Galois theory
(see~\cite{Girstmair.1987, Abdeljaouad.TIATG}), computational
commutative algebra (see~\cite{faugere_rahmany.2009.sagbigroebner,
  Borie_Thiery.2011.Invariants, Borie.2011.Thesis}) and generation of
unlabeled with repetitions species of structures rely on the following
computational building block.

Let $\NN$ be the set of non-negative integers. An \emph{integer vector}
of length $n$ is an element of $\NN^n$. The symmetric group $\sg[n]$
acts on positions on integer vectors in $\NN^n$: for $\sigma$ a
permutation and $(v_1,\dots,v_n)$ an integer vector,
\begin{displaymath}
  \sigma.(v_1,\dots,v_n) := (v_{\sigma^{-1}(1)},\dots,v_{\sigma^{-1}(n)})\,.
\end{displaymath}
This action coincides with the usual action of $\sg[n]$ on monomials
in the multivariate polynomial ring $\KK[\x]$ with $\KK$ a field and
$\x:= x_1, \dots , x_n$ indeterminates.

\begin{problem}\label{enumeration_pb}
  Let $G\subset \sg[n]$ be a permutation group. Enumerate the integer
  vectors of length $n$ modulo the action of $G$.
  
  Note that there are infinitely many such vectors; in practice one
  usually wants to enumerate the vectors with a given sum or content.
\end{problem}

For example, the Problem~\ref{enumeration_pb} contains the listing
non-negative integer matrices with fixed sum up to the permutations of
rows or columns appearing in the theory of multisymmetric
functions~\cite{Gessel87enumerativeapplications, MacMahon} and in the
more recent investigations of multidiagonal coinvariant
\cite{Bergeron.2009.coinvariant.book, MR2820704}.

Define the following equivalence relation over elements of $\NN^n$:
two vectors $\mathbf{u}:= (a_1, \dots , a_n)$ and $\mathbf{v}:= (b_1,
\dots , b_n)$ are equivalent if there exists a permutation $\sigma \in
G$ such that
\begin{displaymath}
  \sigma\cdot \mathbf{u} = (a_{\sigma^{-1}(1)}, \dots , a_{\sigma^{-1}(n)}) =
  (b_1, \dots , b_n) = \mathbf{v}.
\end{displaymath}
Problem~\ref{enumeration_pb} consists in enumerating all $\NN^n /
G$ equivalence classes. 

This problem is not well solved in the literature. Some applications
present a greedy strategy searching and deleting all pairs of vectors
such that the second part can be obtained from the first part. The
most famous sub-problem is the unlabeled graph generation which
consists in enumerate tuples over $0$ and $1$ of length $\binom{n}{2}$
enumerated up to the action of the symmetric groups acting on pair on
nodes. This example has a very efficient implementation in Nauty
which is able to enumerate all graphs over a small number of nodes.

The algorithms presented in this paper have been implemented,
optimized, and intensively tested in \Sage~\cite{Sage}; most features
are integrated in \Sage since release 4.7 (2011-05-26, ticket \#6812, 1303
lines of code including documentation).

\section{Orderly generation and tree structure over integer vectors}
\label{enumeration_strategy}

The orderly strategy consists in setting a total order on objects
before quotienting by the equivalence relation. This allows us to define
a single representative by orbit. Using the lexicographic order on
integer vectors, we will call a vector $\v$ \emph{canonical under the
  action of $G$} or just \emph{canonical} if $\v$ is maximum in its
orbit under $G$ for the lexicographic order:
\begin{displaymath}\label{can_def}
  \v \text{ is canonical } \Leftrightarrow \v = \max_{lex} \{ \sigma\cdot \v \suchthat \sigma \in G \}.
\end{displaymath} 

Now, the goal being to avoid to test systematically if vectors are canonical,
we decided to use a tree structure on the objects in which we will get
properties relaying the \emph{canonical} vectors. Any
result relating fathers, sons and the property of being canonical in
the tree may allowed us to skip some canonical test.

\subsection{Tree Structure over integer vectors}
\label{tree_structure}

Let $\mathbf{r}$ be the vector $\mathbf{r} := (0, \dots , 0)$ called \emph{root},
we build a tree with the following function \emph{father}.

\begin{definition}\label{father} Let $\a = (a_1, a_2, \dots ,
a_n)$ be a tuple of integers of length $n$ which is not the \emph{root}. Let $1 \leqslant i
\leqslant n$ be the position of the last non-zero entry of $\a$. We define the \emph{father} of $\a$
\begin{displaymath}
  father(a_1, a_2, \dots , a_i, 0, 0, \dots , 0) := (a_1, a_2, \dots , a_i - 1, 0, 0, \dots , 0)
\end{displaymath}
\end{definition}

For any integer vector $\mathbf{v} = (a_1, \dots , a_n)$, we can go
back to the generation root $(0, \dots , 0)$ by $sum(\mathbf{v}) :=
a_1 + \dots + a_n$ steps. The corresponding application giving the
children of an integer vector is thus:

\begin{definition}\label{tree_gen_children} Let $\a = (a_1, a_2, \dots ,
a_n)$ be a tuple of integers of length $n$. Let $1 \leqslant i
\leqslant n$ be the position of the last non-zero entry of $\a$ ($i=1$ if
all entries are null). The set of \emph{children} of $\a$ is obtained as:
\begin{displaymath}
  \text{children: }
  (a_1, a_2, \dots , a_i, 0, 0, \dots , 0) \longmapsto
  \left\{
    \begin{array}{c}
    (a_1, a_2, \dots , a_i+1, 0, 0, \dots , 0) \\
    (a_1, a_2, \dots , a_i, 1, 0, \dots ,  0) \\
    (a_1, a_2, \dots , a_i, 0, 1, \dots ,  0) \\
    \dots \\
    (a_1, a_2, \dots , a_i, 0, 0, \dots , 1)
    \end{array}
  \right\}
\end{displaymath}
\end{definition}

\begin{proposition}\label{cut}
  For any permutation group $G \subset \sg[n]$, for any integer vector
  $\mathbf{v}$, if $\mathbf{v}$ is not canonical under $G$, all
  children of $\mathbf{v}$ are not canonical. Therefore, the
  \emph{canonicals} form a "prefix tree" in the tree of all integer
  vectors.

  {\bf Sketch of proof:} When a father is not canonical, there exists
  a permutation such that the permuted vector is greater. Applying the
  same permutation on the children shows also it cannot be canonical.

\end{proposition}

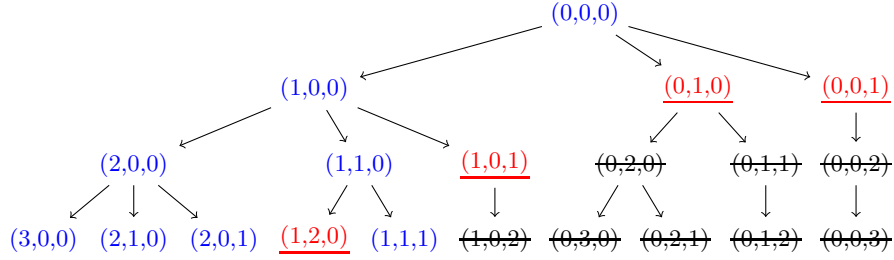
\begin{figure}[h]
  \begin{small}
    \centering
    \begin{tikzpicture}[xscale=1.2]
      \tikzset{
        ed/.style={->, shorten >=1pt},
        cano/.style={color=blue},
        notcano/.style={color=red}
      }

      \draw[cano] (0,0) node (000) {(0,0,0)};
      
      \draw[cano] (-3,-1) node (100) {(1,0,0)};
      \draw[notcano] (1.25,-1) node (010) {\underline{(0,1,0)}};
      \draw[notcano] (3,-1) node (001) {\underline{(0,0,1)}};
      \draw[ed] (000) -- (100);
      \draw[ed] (000) -- (010);
      \draw[ed] (000) -- (001);
      
      \draw[cano] (-5,-2) node (200) {(2,0,0)};
      \draw[cano] (-2.5,-2) node (110) {(1,1,0)};
      \draw[notcano] (-1,-2) node (101) {\underline{(1,0,1)}};
      \draw[ed] (100) -- (200);
      \draw[ed] (100) -- (110);
      \draw[ed] (100) -- (101);
      
      \draw (0.5,-2) node (020) {\sout{(0,2,0)}};
      \draw (2,-2) node (011) {\sout{(0,1,1)}};
      \draw[ed] (010) -- (020);
      \draw[ed] (010) -- (011);
      
      \draw (3,-2) node (002) {\sout{(0,0,2)}};
      \draw[ed] (001) -- (002);
      
      \draw[cano] (-6,-3) node (300) {(3,0,0)};
      \draw[cano] (-5,-3) node (210) {(2,1,0)};
      \draw[cano] (-4,-3) node (201) {(2,0,1)};
      \draw[notcano] (-3,-3) node (120) {\underline{(1,2,0)}};
      \draw[cano] (-2,-3) node (111) {(1,1,1)};
      \draw (-1,-3) node (102) {\sout{(1,0,2)}};
      \draw[ed] (200) -- (300);
      \draw[ed] (200) -- (210);
      \draw[ed] (200) -- (201);
      \draw[ed] (110) -- (120);
      \draw[ed] (110) -- (111);
      \draw[ed] (101) -- (102);

      \draw (0,-3) node (030) {\sout{(0,3,0)}};
      \draw (1,-3) node (021) {\sout{(0,2,1)}};
      \draw[ed] (020) -- (030);
      \draw[ed] (020) -- (021);
      
      \draw (2,-3) node (012) {\sout{(0,1,2)}};
      \draw[ed] (011) -- (012);
      
      \draw (3,-3) node (003) {\sout{(0,0,3)}};
      \draw[ed] (002) -- (003);
      
    \end{tikzpicture}
  \end{small}
  \caption{Enumeration tree of integer vectors modulo the action of
  $G = \langle (1,2,3) \rangle \subset \sg[3]$, the cyclic
  group of degree $3$.}~\label{tree_enumeration}
\end{figure}

Figure~\ref{tree_enumeration} displays integer vectors of length
$3$ whose sum is at most $3$ and shows the tree relations between
them. Choosing the cyclic group of order $3$ and using the generation
strategy, underlined integer vectors are tested but are recognized to
be not \emph{canonical}. Using Proposition~\ref{cut}, crossed-out
integer vectors are not tested as they cannot be \emph{canonical} as
children of non canonical vectors.

Our strategy consists now in making a breath first search over the
sub-tree of \emph{canonicals}. This is done lazily using {\tt Python}
iterators.

\subsection{Testing whether an integer vector is canonical}

As we have seen, the fundamental operation for orderly generation is
to test whether an integer vector is canonical; it is thus vital to
optimize this operation. To this end, we use the work horse of
computational group theory for permutation groups: stabilizer chains
and strong generating sets.

Following the needs required by applications, we want to test
massively if vectors are canonical or not. For this reason, we will
use a \emph{strong generating system} of the group $G$. We can compute
this last item in almost linear
time~\cite{Seress.2003.PermiutationGroupAlgorithms} using \gap~\cite{GAP}.

Let $n$ a positive integer and $G$ a permutation group $G \subset
\sg[n]$. Recall that its \emph{stabilizer chain} is $G_n = \{ e \}
\subset G_{n-1} \subset \dots \subset G_1 \subset G_0 = G$, where
\begin{displaymath}
  \forall i, 1 \leqslant i \leqslant n : G_i := \{ g \in G | \forall j
  \leqslant i : g(j) = j \}\,.
\end{displaymath}

From this chain, we build a \emph{strong generating system} $T =
\{T_1, T_2, \dots , T_n\}$ where $T_i$ is a transversal of $G_{i-1} /
G_i$. This set of strong generators is particularly adapted to the
partial lexicographic order as stabilizers are defined with positions
$1, 2, \dots , n$ from left to right.

Let $n$ and $i$ be two positive integers such that $1 \leqslant i
\leqslant n$. For $\mathbf{v} = (v_1, \dots , v_n)$ and $\mathbf{w} =
(w_1, \dots , w_n)$ two integer vectors of length $n$, let us define the
following binary relations
\begin{displaymath}
  \begin{array}{rcl}
    \mathbf{v} <_i \mathbf{w} & \Longleftrightarrow & (v_1, \dots , v_i) <_{lex} (w_1, \dots , w_i) \\
    \mathbf{v} \leqslant_i \mathbf{w} & \Longleftrightarrow & (v_1, \dots , v_i) \leqslant_{lex} (w_1, \dots , w_i) \\
    \mathbf{v} =_i \mathbf{w} & \Longleftrightarrow & \forall j, 1 \leqslant j \leqslant i : v_j = w_j \\
  \end{array}
\end{displaymath}
where $<_{lex}$ and $\leqslant_{lex}$ represent regular strict and large
lexicographic comparison.

Algorithm~\ref{is_canonical} is a natural extension of McKay's canonical
graph labeling algorithm as it is explained in~\cite{hartke-radcliffe.2009.canonical}.

\begin{algorithm}[H]
\caption{Testing whether an integer vector is canonical}
\label{is_canonical}
  \raggedright

  \renewcommand{\labelitemi}{$\bullet$}

  {\bf Arguments} \\
  $\bullet$ $\mathbf{v}$: An integer vector of length $n$; \\
  $\bullet$ $\sgs$: A strong generating set for $G$, as a list
  $\{T_1,\dots,T_n\}$ of transversals. \\


\begin{displaymath}
    \begin{small}
      \begin{array}{l}
        \textbf{def } is\_canonical(\mathbf{v}, \sgs) :  \\
        \ind todo \leftarrow \{\mathbf{v}\} \\
        \ind \textbf{for } i \in \{1, 2, \dots , n\} : \\ 
        \ind \ind new\_todo \leftarrow \{ \hspace{0.13cm} \} \\
        \ind \ind \textbf{for } \mathbf{w} \in todo : \\
        \ind \ind \ind children \leftarrow \{g \cdot \mathbf{w} | g \in T_i \} \\
        \ind \ind \ind \textbf{for } child \in children : \\
        \ind \ind \ind \ind \textbf{if } \mathbf{v} <_i child : \\
        \ind \ind \ind \ind \ind \textbf{return } False \\
        \ind \ind \ind \ind \textbf{else } : \\
        \ind \ind \ind \ind \ind \textbf{if } \mathbf{v} =_i child \textbf{ and } child \notin new\_todo : \\
        \ind \ind \ind \ind \ind \ind new\_todo \leftarrow new\_todo \cup \{child\} \\
        \ind \ind todo \leftarrow new\_todo \\
        \ind \textbf{return } True \\
      \end{array}
    \end{small}
  \end{displaymath}
\end{algorithm}

Algorithm~\ref{is_canonical} takes advantage of partial lexicographic
orders and the \emph{strong generating system} of the group $G$. It
tries to explore only a small part of the orbit of the vector
$\mathbf{v}$; the worst case complexity of this step is bounded by the
size of the orbit, and not by $|G|$. In this sense, it does take
into account the automorphism group of the vector $\mathbf{v}$.

\begin{proposition}\label{proof_is_canonical}
Let $n$ be a positive integer and $G$ a subgroup of
$\sg[n]$. Let $\mathbf{v}$ be an integer vector of length $n$.
Algorithm~\ref{is_canonical} returns $True$ if $\mathbf{v}$ is
\emph{canonical} under the action of $G$ and returns $False$
otherwise.

  {\bf Sketch of proof:} It is based on the properties of a
  \emph{strong generating system}.

\end{proposition}

\section{Complexity}

\subsection{Theoretical complexity}

\subsubsection{Efficiency of the tree structure}

Let $n$ be a positive integer and $G\subset \sg[n]$ a permutation group. For any non negative integer $d$, let $C(d)$
(resp. $\overline{C}(d)$) be the number of canonical (resp. non
canonical) integer vectors of degree $d$. Based on the tree structure
presented in Section~\ref{tree_structure}, let $T(n)$
(resp. $\overline{T}(n)$) the number of tested (resp. non tested)
integer vectors.

\begin{proposition} \label{compl_NC}
Generating all canonical integer vectors up to degree $d \geqslant 0$ using the
generation strategy presented in Section~\ref{enumeration_strategy}
presents an \emph{absolute error} bounded by $\overline{C}(d)$. Equivalently,
regarding the series, we have
\begin{displaymath}
  \sum_{i=0}^{d} T(i) - \sum_{i=0}^{d} C(i) \leqslant \overline{C}(d)
\end{displaymath}

  {\bf Sketch of proof:} Using Lemma~\ref{cut}, we get this bound
  noticing two tested but non canonical vectors cannot have a
  paternity relation.

\end{proposition}

This \emph{absolute error} is not very explicit (directly usable), but
it can be used to get a \emph{relative error} at the price of a rough
approximation.

\begin{corollary}\label{compl_G}
Let $n$ and $b$ be two positive integers and $G \subset \sg[n]$ a
permutation group. Generating all canonical monomials under the
action of $G$ up to degree $d$ using the generation strategy
presented in Section~\ref{enumeration_strategy} presents a
\emph{relative error} bounded by $\min \{ \frac{n(|G|-1)}{n+d}, n-1
\}$.

  {\bf Sketch of proof:} We use the previous proposition with the fact
  that any integer vector has at least one child but no more than
  $n-1$ children (the generation root is the only one having $n$
  children).

\end{corollary}

The bound is optimal for trivial groups ($\{e\} \subset \sg[n]$), and
seems to be better as the permutation group is of small
cardinality. This relative error becomes better as we go up along the
degree and tends to become optimal when the degree goes to infinity.

\subsubsection{Complexity of testing if a vector is canonical}

We now investigate the complexity of Algorithm~\ref{is_canonical}. We
need first to select a reasonable statistic to collect, which will
define the complexity of this algorithm.

The explosion appearing in the
algorithm is conditioned by the size of the set $new\_todo$. For
$\mathbf{v}$ an integer vector and $\{T_1, \dots , T_n\}$ a
\emph{strong generating system} of a permutation $G$, when $i$
runs over $\{1, 2, \dots , n\}$ in the main loop, the set
$new\_todo_i$ contains at step $i$:
\begin{displaymath}
  new\_todo_i = \{ g_1 \cdots g_i \cdot \mathbf{v} | g_1 \cdots g_i \cdot \mathbf{v} =_i \mathbf{v}, \forall j \leqslant i : g_j \in T_j \}
\end{displaymath}
The right statistic to record is the size of the union of the
$new\_todo_i$ for all $i$ such that the algorithm is still running:
that corresponds to the part of the orbit explored by the
algorithm. This statistic appears to be very difficult to evaluate by a
theoretical way. However, collecting it with a computer is a simple
task.

\subsubsection{Parallelization and memory complexity }

Let us note that this generation engine is trivially amenable for
parallelism: one can devote the study of each branch to a different
processor. Our implementation uses a little framework {\tt
  SearchForest}, co-developed by the author, for exploration trees and
map-reduce operations on them. To get a parallel implementation, it is
sufficient to use the drop-in parallel replacement for {\tt
  SearchForest} under development by Jean-Baptiste Priez and Florent
Hivert.

The memory complexity of the generation engine is reasonable, bounded
by the size of the answer. Indeed, we keep in the cache only the
\emph{Canonical} vectors of degree $d-1$ when we search for those in
degree $d$. In case one wants to only \emph{iterate} through the
elements of a given degree $d$, then this can be achieved with memory
complexity $O(nd)$.

\subsection{Benchmarks design}

To benchmark our implementation, we chose the following problem as test-case.
\begin{problem}\label{canonical_under_staircase}
  Let $n$ be a positive integer and $G\subset \sg[n]$ a permutation
  group. Iterate through all the canonical integer vectors $v$ under
  the staircase (i.e. $v_i \leq n-i$).
\end{problem}
A vector $\mathbf{v}$ of length $n$ is said to be \emph{under the staircase}
when it satisfies $\mathbf{v} \leqslant_{lex} (n-1, n-2 , \dots 1, 0)$.

This problem contains essentially all difficulties that can
appear. The family of $n!$ integer vectors under the staircase
contains vectors with trivial automorphism group as well as vectors
with a lot of symmetries. Applications also require to deal with
this problem as the corresponding family of monomials plays a crucial
role in algebra.

\subsubsection{Benchmarks for transitive permutation groups}

We now need a good family of permutation groups, representative of the
practical use cases. We chose to use the database of all transitive groups
of degree $\leq 30$~\cite{Hulpke.2005.TransitivePermutationGroups} available
in \Sage through the system \gap~\cite{GAP}.

The benchmarks have been run on an off-the-shelf {\tt 2.40 GHz} dual
core Mac Book laptop running {\tt Ubuntu 12.4} and \Sage version 5.3.

\subsection{Benchmarks}

\subsubsection{Tree Structure over integer vectors}

This first benchmark investigates the efficiency of the tree structure
presented in Section~\ref{tree_structure}. As we don't test children
of non canonical integer vectors, one wants to take measures of the
part of tested non canonical vectors (which corresponds to the useless
part of computations). For that, we solve
Problem~\ref{canonical_under_staircase} for each group of the database
and we collect the following information as follows.

\begin{footnotesize}
\begin{table}[H]
  \centering
  \begin{tabular}{|c|c|c|c|c|}
    \hline
    \multicolumn{5}{|c|}{\textbf{Transitive Groups of degree $5$}} \\ \hline
    Database Id. & $|G|$ & Index in $\sg[n]$ & \emph{Canonicals} & number of tests \\ \hline
    1 & 5 & 24 & 71 & 81 \\ \hline
    2 & 10 & 12 & 68 & 81 \\ \hline
    3 & 20 & 6 & 46 & 67 \\ \hline
    4 & 60 & 2 & 41 & 67 \\ \hline
    5 & 120 & 1 & 41 & 67 \\ \hline
  \end{tabular}
  \label{tab_call}
\end{table}
\end{footnotesize}

This table displays the statistics for transitive groups of degree
$5$. \emph{Database Id.} is the integer indexing the group, $|G|$ and
\emph{Index in $\sg[n]$} are respectively the cardinality and the
index of the group $G$ in the symmetric group
$\sg[n]$. \emph{Canonicals} denotes the number of canonical vectors
under the staircase and \emph{number of tests} is the number of times
the algorithm testing if an integer vector is canonical is called.

From this information, we set a quantity $Err$ defined as follows:
\begin{displaymath}
Err := \frac{\text{number of tests} - \text{Canonicals}}{\text{Canonicals}} .
\end{displaymath}

The following figure shows $Err$ depending on the index
$\frac{n!}{|G|}$. The figure contains $166$ crosses, one for each
transitive group over at most $10$ variables. We use a logarithmic
scale on the x axis.

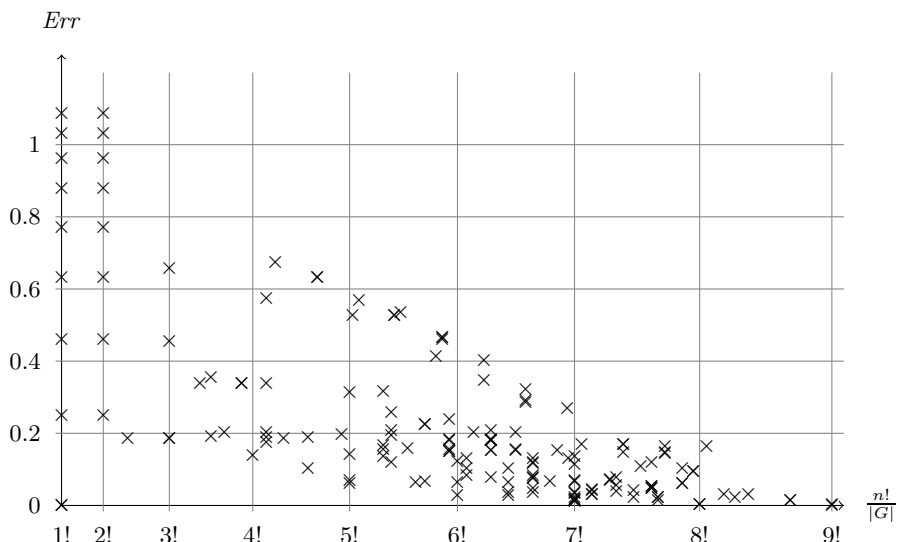
\begin{figure}[H]
  \begin{small}
    \centering
    \begin{tikzpicture}[xscale=0.8, yscale=0.8]
      \tikzset{
        ed/.style={->, shorten >=1pt},
        cano/.style={color=blue},
        notcano/.style={color=red}
      }
    \tikzstyle{vertical}=[very thin, color=gray] 
    \tikzstyle{degree}=[thin] 
    \tikzstyle{horizontal}=[-, color=gray]

      \draw[->] (0,0) -- (13,0); 
      \draw[->] (0,0) -- (0,7.5);

      \draw (13.2, 0) node[right] {$\frac{n!}{|G|}$}; 
      \draw (0, 7.8) node[above] {\emph{Err}};

\draw[horizontal] (-0.1, 6) -- (13, 6);
\draw[horizontal] (-0.1, 4.8) -- (13, 4.8);
\draw[horizontal] (-0.1, 3.6) -- (13, 3.6);
\draw[horizontal] (-0.1, 2.4) -- (13, 2.4);
\draw[horizontal] (-0.1, 1.2) -- (13, 1.2);

\draw (-0.2, 6) node[left] {$1$};
\draw (-0.2, 4.8) node[left] {$0.8$};
\draw (-0.2, 3.6) node[left] {$0.6$};
\draw (-0.2, 2.4) node[left] {$0.4$};
\draw (-0.2, 1.2) node[left] {$0.2$};
\draw (-0.2, 0) node[left] {$0$};

\draw[vertical] (0.693,-0.1) -- (0.693,7.2);
\draw[vertical] (1.791,-0.1) -- (1.791,7.2);
\draw[vertical] (3.178,-0.1) -- (3.178,7.2);
\draw[vertical] (4.787,-0.1) -- (4.787,7.2);
\draw[vertical] (6.579,-0.1) -- (6.579,7.2);
\draw[vertical] (8.525,-0.1) -- (8.525,7.2);
\draw[vertical] (10.604,-0.1) -- (10.604,7.2);
\draw[vertical] (12.801,-0.1) -- (12.801,7.2);

\draw (0, -0.2) node[below] {$1!$};
\draw (0.693, -0.2) node[below] {$2!$};
\draw (1.791, -0.2) node[below] {$3!$};
\draw (3.178, -0.2) node[below] {$4!$};
\draw (4.787, -0.2) node[below] {$5!$};
\draw (6.579, -0.2) node[below] {$6!$};
\draw (8.525, -0.2) node[below] {$7!$};
\draw (10.604, -0.2) node[below] {$8!$};
\draw (12.801, -0.2) node[below] {$9!$};

\draw  (6.73340189183736, 0.787793615739811)  node {$\times$};
\draw  (5.48063892334199, 0.717343173431734)  node {$\times$};
\draw  (2.99573227355399, 2.03112840466926)  node {$\times$};
\draw  (6.04025471127741, 0.404133998574483)  node {$\times$};
\draw  (7.42654907239730, 0.173674426838984)  node {$\times$};
\draw  (9.80609520652748, 0.292063075578712)  node {$\times$};
\draw  (4.65396035015752, 1.19510774606872)  node {$\times$};
\draw  (9.21830854162536, 0.345764910522126)  node {$\times$};
\draw  (7.83201418050547, 0.794068823186325)  node {$\times$};
\draw  (7.54433210805369, 1.21361786219915)  node {$\times$};
\draw  (10.0292387578417, 0.893734820510868)  node {$\times$};
\draw  (6.32793678372919, 2.77106153790167)  node {$\times$};
\draw  (7.13886699994552, 0.920903110047847)  node {$\times$};
\draw  (4.78749174278205, 0.421875000000000)  node {$\times$};
\draw  (0.000000000000000, 5.78166550034990)  node {$\times$};
\draw  (5.34710753071747, 1.00356188780053)  node {$\times$};
\draw  (8.11969625295725, 0.401679301910032)  node {$\times$};
\draw  (8.52516136106541, 0.111299017017229)  node {$\times$};
\draw  (5.48063892334199, 1.25638911788953)  node {$\times$};
\draw  (10.7223859384016, 0.991922514580619)  node {$\times$};
\draw  (6.73340189183736, 0.504701553556827)  node {$\times$};
\draw  (3.40119738166216, 1.22891566265060)  node {$\times$};
\draw  (0.000000000000000, 3.80487804878049)  node {$\times$};
\draw  (11.0100680108534, 0.195809470351433)  node {$\times$};
\draw  (7.71423114484909, 1.72374519931135)  node {$\times$};
\draw  (0.693147180559945, 5.27102803738318)  node {$\times$};
\draw  (2.99573227355399, 2.03112840466926)  node {$\times$};
\draw  (7.83201418050547, 0.466635804510284)  node {$\times$};
\draw  (9.80609520652748, 0.289875820132713)  node {$\times$};
\draw  (5.34710753071747, 1.90166253979484)  node {$\times$};
\draw  (8.81284343351719, 0.194435322566200)  node {$\times$};
\draw  (6.44571981938558, 1.44429400386847)  node {$\times$};
\draw  (5.52942908751142, 3.16625041878598)  node {$\times$};
\draw  (7.02108396428914, 2.07938240798503)  node {$\times$};
\draw  (7.13886699994552, 1.08404099560761)  node {$\times$};
\draw  (0.693147180559945, 4.62595419847328)  node {$\times$};
\draw  (0.000000000000000, 0.000000000000000)  node {$\times$};
\draw  (7.13886699994552, 0.476255088195387)  node {$\times$};
\draw  (8.52516136106541, 0.412933909462634)  node {$\times$};
\draw  (4.24849524204936, 3.80355837192298)  node {$\times$};
\draw  (0.000000000000000, 1.50000000000000)  node {$\times$};
\draw  (8.52516136106541, 0.227000825082508)  node {$\times$};
\draw  (8.52516136106541, 0.0856412441414572)  node {$\times$};
\draw  (0.000000000000000, 2.76923076923077)  node {$\times$};
\draw  (7.54433210805369, 0.936239113658913)  node {$\times$};
\draw  (10.3169208302935, 0.617293110105774)  node {$\times$};
\draw  (9.91145572218531, 0.134477078227329)  node {$\times$};
\draw  (7.13886699994552, 1.08788740829465)  node {$\times$};
\draw  (0.693147180559945, 6.52944328669247)  node {$\times$};
\draw  (9.33609157728174, 1.02678099147325)  node {$\times$};
\draw  (8.52516136106541, 0.693750000000000)  node {$\times$};
\draw  (3.40119738166216, 2.03311258278146)  node {$\times$};
\draw  (6.22257626807137, 2.48054610393727)  node {$\times$};
\draw  (7.83201418050547, 0.450851035756696)  node {$\times$};
\draw  (4.83628190695148, 3.16625041878598)  node {$\times$};
\draw  (9.33609157728174, 0.893734820510868)  node {$\times$};
\draw  (1.79175946922805, 1.12500000000000)  node {$\times$};
\draw  (7.54433210805369, 0.925776738359699)  node {$\times$};
\draw  (9.50599061407714, 0.249716645855345)  node {$\times$};
\draw  (6.44571981938558, 0.912342421577250)  node {$\times$};
\draw  (3.68887945411394, 1.11724137931034)  node {$\times$};
\draw  (1.79175946922805, 2.73913043478261)  node {$\times$};
\draw  (5.48063892334199, 1.54742466973430)  node {$\times$};
\draw  (11.4155331189616, 0.186881947341300)  node {$\times$};
\draw  (6.57925121201010, 0.737188721959251)  node {$\times$};
\draw  (3.40119738166216, 1.05882352941176)  node {$\times$};
\draw  (7.42654907239730, 0.627806385169928)  node {$\times$};
\draw  (5.48063892334199, 1.16801470588235)  node {$\times$};
\draw  (4.78749174278205, 1.88346186803770)  node {$\times$};
\draw  (0.693147180559945, 1.50000000000000)  node {$\times$};
\draw  (9.21830854162536, 0.237348006249575)  node {$\times$};
\draw  (6.44571981938558, 1.10113713382032)  node {$\times$};
\draw  (10.0292387578417, 0.991922514580619)  node {$\times$};
\draw  (6.32793678372919, 2.79731357447896)  node {$\times$};
\draw  (7.13886699994552, 1.24658162797100)  node {$\times$};
\draw  (2.30258509299405, 2.03112840466926)  node {$\times$};
\draw  (10.0292387578417, 0.877519988524058)  node {$\times$};
\draw  (0.693147180559945, 5.78166550034990)  node {$\times$};
\draw  (8.64294439672180, 1.02678099147325)  node {$\times$};
\draw  (7.83201418050547, 0.723441326050104)  node {$\times$};
\draw  (0.000000000000000, 4.62595419847328)  node {$\times$};
\draw  (8.81284343351719, 0.279547885948184)  node {$\times$};
\draw  (8.52516136106541, 0.108336413072171)  node {$\times$};
\draw  (9.21830854162536, 0.470212934629227)  node {$\times$};
\draw  (5.75257263882563, 0.950378304018866)  node {$\times$};
\draw  (10.3169208302935, 0.365645782058779)  node {$\times$};
\draw  (10.6046029027453, 0.0288751924038630)  node {$\times$};
\draw  (6.73340189183736, 0.615771428571429)  node {$\times$};
\draw  (3.40119738166216, 1.13414634146341)  node {$\times$};
\draw  (0.000000000000000, 6.19255297263937)  node {$\times$};
\draw  (8.52516136106541, 0.817105592579788)  node {$\times$};
\draw  (11.1923895676474, 0.144032353681612)  node {$\times$};
\draw  (7.71423114484909, 1.76083922590249)  node {$\times$};
\draw  (0.000000000000000, 5.27102803738318)  node {$\times$};
\draw  (0.000000000000000, 0.000000000000000)  node {$\times$};
\draw  (7.83201418050547, 0.703680000000000)  node {$\times$};
\draw  (12.1086802995215, 0.0843593016110388)  node {$\times$};
\draw  (5.34710753071747, 0.944506104328524)  node {$\times$};
\draw  (9.91145572218531, 0.0817636674010417)  node {$\times$};
\draw  (6.44571981938558, 0.936311323515069)  node {$\times$};
\draw  (3.17805383034795, 0.845070422535211)  node {$\times$};
\draw  (12.1086802995215, 0.0882134302469646)  node {$\times$};
\draw  (7.02108396428914, 2.42075094995814)  node {$\times$};
\draw  (5.88610403145016, 0.394569367840475)  node {$\times$};
\draw  (1.79175946922805, 3.94285714285714)  node {$\times$};
\draw  (7.42654907239730, 0.216391716391716)  node {$\times$};
\draw  (8.52516136106541, 0.421776916689731)  node {$\times$};
\draw  (4.24849524204936, 3.79639551875304)  node {$\times$};
\draw  (8.52516136106541, 0.200783741911966)  node {$\times$};
\draw  (8.52516136106541, 0.0718189768255356)  node {$\times$};
\draw  (0.000000000000000, 6.52944328669247)  node {$\times$};
\draw  (0.693147180559945, 2.76923076923077)  node {$\times$};
\draw  (8.41980084540759, 0.789391427894861)  node {$\times$};
\draw  (10.3169208302935, 0.365267416202585)  node {$\times$};
\draw  (6.32793678372919, 2.80523941006128)  node {$\times$};
\draw  (7.13886699994552, 1.08431616731375)  node {$\times$};
\draw  (4.09434456222210, 0.629032258064516)  node {$\times$};
\draw  (6.04025471127741, 1.35570263198399)  node {$\times$};
\draw  (9.11294802596753, 0.431407454423290)  node {$\times$};
\draw  (8.40737832540903, 1.62493690055527)  node {$\times$};
\draw  (7.83201418050547, 0.218305468971380)  node {$\times$};
\draw  (9.62377364973352, 0.653596044232099)  node {$\times$};
\draw  (1.79175946922805, 1.12500000000000)  node {$\times$};
\draw  (10.4992423870874, 0.573247566145644)  node {$\times$};
\draw  (9.50599061407714, 0.141494072275389)  node {$\times$};
\draw  (6.57925121201010, 0.382457526669301)  node {$\times$};
\draw  (4.09434456222210, 1.13414634146341)  node {$\times$};
\draw  (0.693147180559945, 3.80487804878049)  node {$\times$};
\draw  (8.23747928861363, 0.925776738359699)  node {$\times$};
\draw  (5.34710753071747, 0.823754789272031)  node {$\times$};
\draw  (7.71423114484909, 1.94032439527470)  node {$\times$};
\draw  (6.57925121201010, 0.177399756986634)  node {$\times$};
\draw  (2.70805020110221, 1.22891566265060)  node {$\times$};
\draw  (7.83201418050547, 0.499586435070306)  node {$\times$};
\draw  (9.80609520652748, 0.324626138294021)  node {$\times$};
\draw  (4.94164242260930, 3.42035087719298)  node {$\times$};
\draw  (8.81284343351719, 0.195302506599897)  node {$\times$};
\draw  (6.44571981938558, 0.898834951456311)  node {$\times$};
\draw  (6.85118492749374, 1.21361786219915)  node {$\times$};
\draw  (12.8018274800815, 0.0112082427887992)  node {$\times$};
\draw  (5.63478960316925, 3.21717847718842)  node {$\times$};
\draw  (7.13886699994552, 1.10944860072042)  node {$\times$};
\draw  (2.48490664978800, 2.13658536585366)  node {$\times$};
\draw  (3.40119738166216, 3.45418326693227)  node {$\times$};
\draw  (9.11294802596753, 0.428258672123468)  node {$\times$};
\draw  (3.55534806148941, 4.04595404595405)  node {$\times$};
\draw  (0.693147180559945, 6.19255297263937)  node {$\times$};
\draw  (8.81284343351719, 0.259960206923995)  node {$\times$};
\draw  (8.52516136106541, 0.137949260042283)  node {$\times$};
\draw  (6.04025471127741, 1.34613022113022)  node {$\times$};
\draw  (10.4992423870874, 0.571792079842086)  node {$\times$};
\draw  (10.6046029027453, 0.0269320843091335)  node {$\times$};
\draw  (7.13886699994552, 0.926441351888668)  node {$\times$};
\draw  (4.78749174278205, 0.371134020618557)  node {$\times$};
\draw  (9.80609520652748, 0.715021971633777)  node {$\times$};
\draw  (9.33609157728174, 1.02045014042166)  node {$\times$};
\draw  (7.42654907239730, 0.384227418615314)  node {$\times$};
\draw  (4.78749174278205, 0.855036855036855)  node {$\times$};
\draw  (7.83201418050547, 0.287544483985765)  node {$\times$};
\draw  (5.52942908751142, 3.16625041878598)  node {$\times$};
\draw  (9.80609520652748, 0.318541372405800)  node {$\times$};
\draw  (1.09861228866811, 1.12500000000000)  node {$\times$};
\draw  (9.91145572218531, 0.143394963607463)  node {$\times$};
\draw  (6.44571981938558, 1.08431616731375)  node {$\times$};
\draw  (2.48490664978800, 1.14705882352941)  node {$\times$};
\draw  (12.8018274800815, 0.0198779805632026)  node {$\times$};

    \end{tikzpicture}
  \end{small}
  \caption{Relative Error between number of tested vectors and number
    of \emph{canonicals} vectors.}~\label{call_pic}
  \label{call_pic}
\end{figure}

\subsubsection{Empirical complexity of testing if a vector is canonical}

Algorithm~\ref{is_canonical} needs to explore a part of the orbit of the
tested integer vectors. The following table displays for each
transitive group over $5$ variables, the number of elements of all
orbits of tested vectors solving
Problem~\ref{canonical_under_staircase} compared to the total number
of integer vectors explored.

\begin{table}[H]
  \centering
  \begin{tabular}{|c|c|c|c|c|}
    \hline
    \multicolumn{5}{|c|}{\textbf{Transitive Groups of degree $5$}} \\ \hline
    Database Id. & $|G|$ & Index in $\sg[n]$ & total orbits & total explored \\ \hline
    1 & 5 & 24 & 401 & 351 \\ \hline
    2 & 10 & 12 & 691 & 393 \\ \hline
    3 & 20 & 6 & 1091 & 365 \\ \hline
    4 & 60 & 2 & 1891 & 328 \\ \hline
    5 & 120 & 1 & 1891 & 326 \\ \hline
  \end{tabular}
  \label{tab_call}
\end{table}

Now we define $Ratio$ to be the average size of the orbit needed to be
explored to know if an integer vector is canonical:
\begin{displaymath}
  Ratio := \frac{\text{total explored}}{\text{total orbits}} .
\end{displaymath}
The following figure plots $Ratio$ in terms of $|G|$ for transitive
groups on at most $9$ variables.

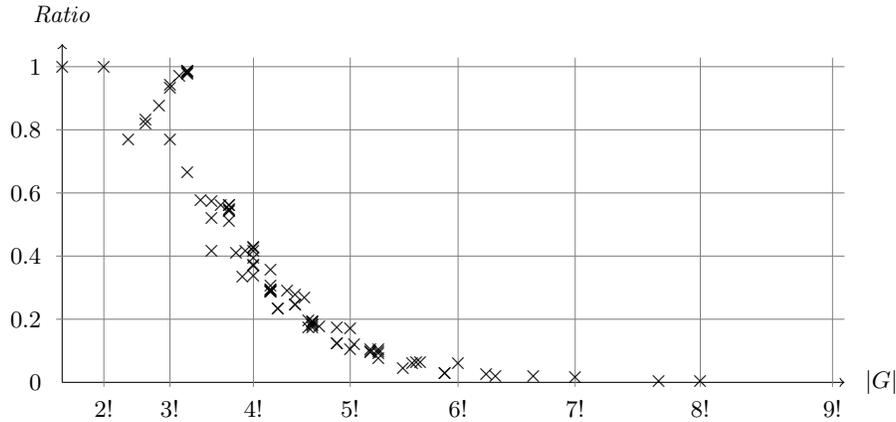
\begin{figure}[H]
  \begin{small}
    \centering
    \begin{tikzpicture}[xscale=0.8, yscale=0.6]
      \tikzset{
        ed/.style={->, shorten >=1pt},
        cano/.style={color=blue},
        notcano/.style={color=red}
      }
    \tikzstyle{vertical}=[very thin, color=gray] 
    \tikzstyle{degree}=[thin] 
    \tikzstyle{horizontal}=[-, color=gray]

      \draw[->] (0,0) -- (13,0); 
      \draw[->] (0,0) -- (0,7.5);

      \draw (13.2, 0) node[right] {$|G|$}; 
      \draw (0, 7.8) node[above] {\emph{Ratio}};

\draw[horizontal] (-0.1, 7) -- (13, 7);
\draw[horizontal] (-0.1, 5.6) -- (13, 5.6);
\draw[horizontal] (-0.1, 4.2) -- (13, 4.2);
\draw[horizontal] (-0.1, 2.8) -- (13, 2.8);
\draw[horizontal] (-0.1, 1.4) -- (13, 1.4);

\draw (-0.2, 7) node[left] {$1$};
\draw (-0.2, 5.6) node[left] {$0.8$};
\draw (-0.2, 4.2) node[left] {$0.6$};
\draw (-0.2, 2.8) node[left] {$0.4$};
\draw (-0.2, 1.4) node[left] {$0.2$};
\draw (-0.2, 0) node[left] {$0$};

\draw[vertical] (0.693,-0.1) -- (0.693,7.2);
\draw[vertical] (1.791,-0.1) -- (1.791,7.2);
\draw[vertical] (3.178,-0.1) -- (3.178,7.2);
\draw[vertical] (4.787,-0.1) -- (4.787,7.2);
\draw[vertical] (6.579,-0.1) -- (6.579,7.2);
\draw[vertical] (8.525,-0.1) -- (8.525,7.2);
\draw[vertical] (10.604,-0.1) -- (10.604,7.2);
\draw[vertical] (12.801,-0.1) -- (12.801,7.2);

\draw (0.693, -0.2) node[below] {$2!$};
\draw (1.791, -0.2) node[below] {$3!$};
\draw (3.178, -0.2) node[below] {$4!$};
\draw (4.787, -0.2) node[below] {$5!$};
\draw (6.579, -0.2) node[below] {$6!$};
\draw (8.525, -0.2) node[below] {$7!$};
\draw (10.604, -0.2) node[below] {$8!$};
\draw (12.801, -0.2) node[below] {$9!$};

\draw  (3.04452243772342, 2.90922112232231)  node {$\times$};
\draw  (3.58351893845611, 1.63997113997114)  node {$\times$};
\draw  (5.95064255258773, 0.449683429752578)  node {$\times$};
\draw  (10.6046029027453, 0.0308688919009269)  node {$\times$};
\draw  (2.07944154167984, 6.87241914232120)  node {$\times$};
\draw  (3.87120101090789, 1.72299273021884)  node {$\times$};
\draw  (3.17805383034795, 2.60789106145251)  node {$\times$};
\draw  (4.78749174278205, 1.20676890534109)  node {$\times$};
\draw  (7.83201418050547, 0.129027941035426)  node {$\times$};
\draw  (3.58351893845611, 1.64332313965341)  node {$\times$};
\draw  (5.25749537202778, 0.544451170242343)  node {$\times$};
\draw  (0.000000000000000, 7.00000000000000)  node {$\times$};
\draw  (3.46573590279973, 2.49037646118751)  node {$\times$};
\draw  (1.79175946922805, 5.38461538461539)  node {$\times$};
\draw  (2.07944154167984, 6.90136892674303)  node {$\times$};
\draw  (3.17805383034795, 2.92052980132450)  node {$\times$};
\draw  (3.46573590279973, 2.04463515675492)  node {$\times$};
\draw  (5.12396397940326, 0.727142455745996)  node {$\times$};
\draw  (5.25749537202778, 0.644533662472220)  node {$\times$};
\draw  (3.17805383034795, 2.76370765672119)  node {$\times$};
\draw  (1.09861228866811, 5.38461538461539)  node {$\times$};
\draw  (4.15888308335967, 1.28631954325334)  node {$\times$};
\draw  (3.46573590279973, 2.05061318499618)  node {$\times$};
\draw  (4.27666611901606, 1.24987869966036)  node {$\times$};
\draw  (6.57925121201010, 0.426700121766724)  node {$\times$};
\draw  (4.85203026391962, 0.851176759234973)  node {$\times$};
\draw  (0.693147180559945, 7.00000000000000)  node {$\times$};
\draw  (2.77258872223978, 3.58874328131415)  node {$\times$};
\draw  (4.15888308335967, 1.34609700550565)  node {$\times$};
\draw  (1.60943791243410, 6.12718204488778)  node {$\times$};
\draw  (2.63905732961526, 3.94065299215126)  node {$\times$};
\draw  (4.78749174278205, 0.740480108960899)  node {$\times$};
\draw  (6.35610766069589, 0.210688195147544)  node {$\times$};
\draw  (4.56434819146784, 1.21318176964441)  node {$\times$};
\draw  (2.77258872223978, 3.92888749183935)  node {$\times$};
\draw  (1.38629436111989, 5.81690140845070)  node {$\times$};
\draw  (4.02535169073515, 1.88885534747548)  node {$\times$};
\draw  (2.48490664978800, 3.65099292398996)  node {$\times$};
\draw  (4.09434456222210, 1.21417239555791)  node {$\times$};
\draw  (1.94591014905531, 6.79890937710816)  node {$\times$};
\draw  (3.87120101090789, 1.71893491124260)  node {$\times$};
\draw  (5.66296048013595, 0.306634834543845)  node {$\times$};
\draw  (7.20340552108309, 0.148764505495203)  node {$\times$};
\draw  (2.07944154167984, 6.84224534863389)  node {$\times$};
\draw  (3.46573590279973, 2.15154131661545)  node {$\times$};
\draw  (1.79175946922805, 6.58873929008568)  node {$\times$};
\draw  (3.73766961828337, 2.04234302809297)  node {$\times$};
\draw  (5.25749537202778, 0.733366962281547)  node {$\times$};
\draw  (3.17805383034795, 2.99964633068081)  node {$\times$};
\draw  (3.46573590279973, 2.01464036184997)  node {$\times$};
\draw  (1.79175946922805, 6.53101736972705)  node {$\times$};
\draw  (5.12396397940326, 0.667401502060417)  node {$\times$};
\draw  (2.77258872223978, 3.83087403811678)  node {$\times$};
\draw  (4.15888308335967, 1.21424177768898)  node {$\times$};
\draw  (3.46573590279973, 2.00685357354578)  node {$\times$};
\draw  (5.88610403145016, 0.440647946570237)  node {$\times$};
\draw  (6.35610766069589, 0.210245110809868)  node {$\times$};
\draw  (4.56434819146784, 0.860485576915980)  node {$\times$};
\draw  (2.77258872223978, 3.84057605944795)  node {$\times$};
\draw  (1.38629436111989, 5.74626865671642)  node {$\times$};
\draw  (4.15888308335967, 1.34735358173192)  node {$\times$};
\draw  (2.89037175789616, 2.87101996762008)  node {$\times$};
\draw  (2.99573227355399, 2.34188817598533)  node {$\times$};
\draw  (3.17805383034795, 2.58249211356467)  node {$\times$};
\draw  (5.81711115996320, 0.421628173760619)  node {$\times$};
\draw  (9.91145572218531, 0.0362919723042597)  node {$\times$};
\draw  (2.07944154167984, 6.87559031589614)  node {$\times$};
\draw  (3.87120101090789, 1.95610556235886)  node {$\times$};
\draw  (3.17805383034795, 2.37207935341661)  node {$\times$};
\draw  (8.52516136106541, 0.121065180879562)  node {$\times$};
\draw  (5.25749537202778, 0.686318902185464)  node {$\times$};
\draw  (3.17805383034795, 2.98651743012104)  node {$\times$};
\draw  (2.07944154167984, 6.91716144197451)  node {$\times$};
\draw  (2.48490664978800, 2.92052980132450)  node {$\times$};
\draw  (3.46573590279973, 2.00638933388455)  node {$\times$};
\draw  (2.48490664978800, 4.01697530864197)  node {$\times$};
\draw  (5.12396397940326, 0.685192643180614)  node {$\times$};
\draw  (2.77258872223978, 3.79453479847479)  node {$\times$};
\draw  (4.15888308335967, 1.33687839299615)  node {$\times$};
\draw  (3.46573590279973, 2.06211401821641)  node {$\times$};
\draw  (4.09434456222210, 1.38393482831114)  node {$\times$};
\draw  (7.04925484125584, 0.174218249430101)  node {$\times$};
\draw  (4.56434819146784, 0.863787675112280)  node {$\times$};
\draw  (2.77258872223978, 3.92344193817878)  node {$\times$};
\draw  (2.07944154167984, 4.66666666666667)  node {$\times$};
\draw  (4.15888308335967, 1.27147470025307)  node {$\times$};
\draw  (2.30258509299405, 4.03964757709251)  node {$\times$};

    \end{tikzpicture}
  \end{small}
  \caption{Average, over all integer vectors $v$ under the stair case,
    of the number of vectors in the orbit of $v$ explored by
    \texttt{is\_canonical(v)}.}~\label{can_pic}
\end{figure}

\subsubsection{Overall empirical complexity of the generation engine}

We now evaluate the overall complexity by comparing the ratio between
the computations and the size of the output. We define the measure
\emph{Complexity} as follows:
\begin{displaymath}
  Complexity := \frac{\text{total explored}}{\text{Canonicals}} .
\end{displaymath}
The following graph displays \emph{Complexity} in terms of the
size of the group $|G|$ for transitive Groups on up to $9$ variables
(and excluding the alternate and symmetric group of degree $9$).

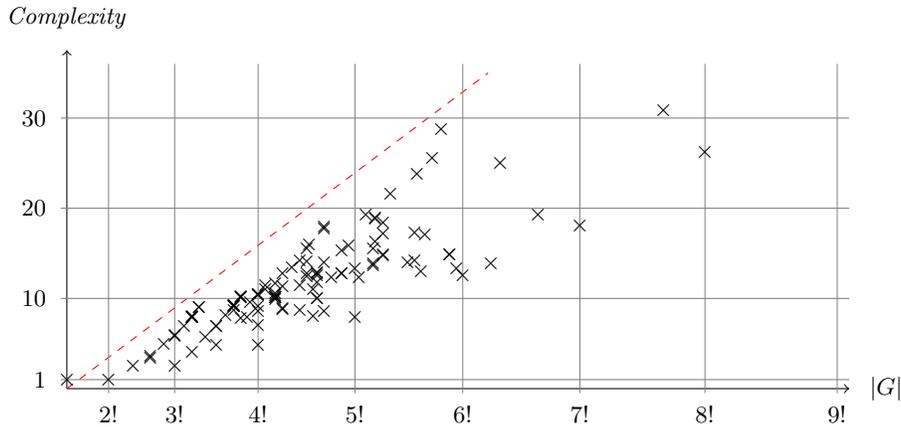
\begin{figure}[H]
  \begin{small}
    \centering
    \begin{tikzpicture}[xscale=0.8, yscale=0.6]
      \tikzset{
        ed/.style={->, shorten >=1pt},
        cano/.style={color=blue},
        notcano/.style={color=red}
      }
    \tikzstyle{vertical}=[very thin, color=gray] 
    \tikzstyle{degree}=[thin] 
    \tikzstyle{horizontal}=[-, color=gray]

      \draw[->] (0,0) -- (13,0); 
      \draw[->] (0,0) -- (0,7.5);

      \draw (13.2, 0) node[right] {$|G|$}; 
      \draw (0, 7.8) node[above] {\emph{Complexity}};

      \draw[dashed, color=red] (0, 0) -- (7.0, 7.0);

\draw[horizontal] (-0.1, 6) -- (13, 6);
\draw[horizontal] (-0.1, 4) -- (13, 4);
\draw[horizontal] (-0.1, 2) -- (13, 2);
\draw[horizontal] (-0.1, 0.2) -- (13, 0.2);

\draw (-0.2, 6) node[left] {$30$};
\draw (-0.2, 4) node[left] {$20$};
\draw (-0.2, 2) node[left] {$10$};
\draw (-0.2, 0.2) node[left] {$1$};

\draw[vertical] (0.693,-0.1) -- (0.693,7.2);
\draw[vertical] (1.791,-0.1) -- (1.791,7.2);
\draw[vertical] (3.178,-0.1) -- (3.178,7.2);
\draw[vertical] (4.787,-0.1) -- (4.787,7.2);
\draw[vertical] (6.579,-0.1) -- (6.579,7.2);
\draw[vertical] (8.525,-0.1) -- (8.525,7.2);
\draw[vertical] (10.604,-0.1) -- (10.604,7.2);
\draw[vertical] (12.801,-0.1) -- (12.801,7.2);

\draw (0.693, -0.2) node[below] {$2!$};
\draw (1.791, -0.2) node[below] {$3!$};
\draw (3.178, -0.2) node[below] {$4!$};
\draw (4.787, -0.2) node[below] {$5!$};
\draw (6.579, -0.2) node[below] {$6!$};
\draw (8.525, -0.2) node[below] {$7!$};
\draw (10.604, -0.2) node[below] {$8!$};
\draw (12.801, -0.2) node[below] {$9!$};

\draw  (6.06842558824411, 5.11430435655491)  node {$\times$};
\draw  (3.04452243772342, 1.92723247232472)  node {$\times$};
\draw  (3.58351893845611, 1.76887159533074)  node {$\times$};
\draw  (5.95064255258773, 3.42465540671714)  node {$\times$};
\draw  (10.6046029027453, 5.24520643806858)  node {$\times$};
\draw  (2.07944154167984, 1.59876334425536)  node {$\times$};
\draw  (3.87120101090789, 2.30404742436631)  node {$\times$};
\draw  (3.17805383034795, 1.71405622489960)  node {$\times$};
\draw  (4.78749174278205, 1.59024390243902)  node {$\times$};
\draw  (5.08759633523238, 2.78693771244427)  node {$\times$};
\draw  (7.83201418050547, 3.86168224299065)  node {$\times$};
\draw  (3.58351893845611, 1.79221789883268)  node {$\times$};
\draw  (5.25749537202778, 2.95889635656173)  node {$\times$};
\draw  (0.000000000000000, 0.200000000000000)  node {$\times$};
\draw  (3.46573590279973, 2.34094301221167)  node {$\times$};
\draw  (1.79175946922805, 0.500000000000000)  node {$\times$};
\draw  (4.27666611901606, 3.55968784378438)  node {$\times$};
\draw  (2.07944154167984, 1.59798679164891)  node {$\times$};
\draw  (3.17805383034795, 0.969230769230769)  node {$\times$};
\draw  (2.89037175789616, 2.04324668193446)  node {$\times$};
\draw  (3.46573590279973, 2.06381194789639)  node {$\times$};
\draw  (4.27666611901606, 2.80093325791855)  node {$\times$};
\draw  (5.12396397940326, 3.26445916114790)  node {$\times$};
\draw  (5.25749537202778, 2.98708815672306)  node {$\times$};
\draw  (3.17805383034795, 2.08361139718503)  node {$\times$};
\draw  (1.09861228866811, 0.500000000000000)  node {$\times$};
\draw  (3.58351893845611, 2.54909539659896)  node {$\times$};
\draw  (4.15888308335967, 2.00546727916259)  node {$\times$};
\draw  (3.46573590279973, 2.03556942277691)  node {$\times$};
\draw  (4.27666611901606, 1.71828793774319)  node {$\times$};
\draw  (4.96981329957600, 3.85597099265398)  node {$\times$};
\draw  (6.57925121201010, 2.52213740458015)  node {$\times$};
\draw  (4.85203026391962, 2.46901837476663)  node {$\times$};
\draw  (0.693147180559945, 0.200000000000000)  node {$\times$};
\draw  (2.77258872223978, 1.73825600000000)  node {$\times$};
\draw  (2.89037175789616, 2.02908854812588)  node {$\times$};
\draw  (4.15888308335967, 2.50519847849410)  node {$\times$};
\draw  (1.60943791243410, 0.988732394366197)  node {$\times$};
\draw  (5.78074351579233, 2.84311843884846)  node {$\times$};
\draw  (2.63905732961526, 1.63716588884175)  node {$\times$};
\draw  (4.78749174278205, 2.66285714285714)  node {$\times$};
\draw  (6.35610766069589, 2.98436434486118)  node {$\times$};
\draw  (4.56434819146784, 3.06560940841055)  node {$\times$};
\draw  (2.77258872223978, 1.84173420232360)  node {$\times$};
\draw  (1.38629436111989, 0.737500000000000)  node {$\times$};
\draw  (3.29583686600433, 2.20934152263962)  node {$\times$};
\draw  (4.02535169073515, 3.20208086395364)  node {$\times$};
\draw  (2.48490664978800, 1.39329268292683)  node {$\times$};
\draw  (4.09434456222210, 1.60000000000000)  node {$\times$};
\draw  (5.08759633523238, 3.10039377461091)  node {$\times$};
\draw  (1.94591014905531, 1.39951397326853)  node {$\times$};
\draw  (3.87120101090789, 1.73333333333333)  node {$\times$};
\draw  (5.66296048013595, 2.80439766081871)  node {$\times$};
\draw  (7.20340552108309, 4.99890438247012)  node {$\times$};
\draw  (3.98898404656427, 2.49841305212893)  node {$\times$};
\draw  (2.07944154167984, 1.59552854122622)  node {$\times$};
\draw  (2.19722457733622, 1.79988030184752)  node {$\times$};
\draw  (3.46573590279973, 1.99098740888005)  node {$\times$};
\draw  (1.79175946922805, 1.18917525773196)  node {$\times$};
\draw  (5.37527840768417, 4.33014672168730)  node {$\times$};
\draw  (3.73766961828337, 2.69975429975430)  node {$\times$};
\draw  (5.25749537202778, 3.43942103022563)  node {$\times$};
\draw  (3.17805383034795, 2.10216871735859)  node {$\times$};
\draw  (3.58351893845611, 2.26989985946913)  node {$\times$};
\draw  (3.46573590279973, 2.03314892344498)  node {$\times$};
\draw  (1.79175946922805, 1.17500000000000)  node {$\times$};
\draw  (4.68213122712422, 3.18539978835245)  node {$\times$};
\draw  (5.12396397940326, 3.77761747732894)  node {$\times$};
\draw  (2.77258872223978, 1.83712717412869)  node {$\times$};
\draw  (3.98898404656427, 2.54779831029644)  node {$\times$};
\draw  (4.15888308335967, 2.35092843326886)  node {$\times$};
\draw  (5.78074351579233, 3.46650031191516)  node {$\times$};
\draw  (3.46573590279973, 2.03215226939971)  node {$\times$};
\draw  (5.88610403145016, 2.60458015267176)  node {$\times$};
\draw  (6.35610766069589, 2.98406044357787)  node {$\times$};
\draw  (4.56434819146784, 2.56342927505002)  node {$\times$};
\draw  (2.77258872223978, 1.83925319963056)  node {$\times$};
\draw  (1.38629436111989, 0.687500000000000)  node {$\times$};
\draw  (3.29583686600433, 2.30123907405629)  node {$\times$};
\draw  (4.15888308335967, 2.58581061272354)  node {$\times$};
\draw  (2.89037175789616, 1.57241379310345)  node {$\times$};
\draw  (2.99573227355399, 1.58695652173913)  node {$\times$};
\draw  (6.22257626807137, 5.75142210331344)  node {$\times$};
\draw  (3.17805383034795, 1.83450980392157)  node {$\times$};
\draw  (5.81711115996320, 4.76612396458155)  node {$\times$};
\draw  (9.91145572218531, 6.16668999300210)  node {$\times$};
\draw  (3.98898404656427, 2.81214009953951)  node {$\times$};
\draw  (2.07944154167984, 1.59661052742727)  node {$\times$};
\draw  (2.19722457733622, 1.79995843427094)  node {$\times$};
\draw  (3.87120101090789, 2.84262857142857)  node {$\times$};
\draw  (3.17805383034795, 1.40609756097561)  node {$\times$};
\draw  (5.08759633523238, 2.73979344994418)  node {$\times$};
\draw  (8.52516136106541, 3.62336448598131)  node {$\times$};
\draw  (5.25749537202778, 3.68975212726600)  node {$\times$};
\draw  (3.17805383034795, 2.10257985257985)  node {$\times$};
\draw  (4.27666611901606, 3.58793401986695)  node {$\times$};
\draw  (2.07944154167984, 1.59920201136861)  node {$\times$};
\draw  (2.48490664978800, 0.969230769230769)  node {$\times$};
\draw  (6.47389069635227, 2.66911576537141)  node {$\times$};
\draw  (3.46573590279973, 2.03154347666838)  node {$\times$};
\draw  (2.48490664978800, 1.39946236559140)  node {$\times$};
\draw  (4.39444915467244, 2.46198889449773)  node {$\times$};
\draw  (5.12396397940326, 3.79768382352941)  node {$\times$};
\draw  (2.77258872223978, 1.82457954232148)  node {$\times$};
\draw  (3.98898404656427, 3.10924802003103)  node {$\times$};
\draw  (4.15888308335967, 2.56776699029126)  node {$\times$};
\draw  (3.46573590279973, 2.11916505033712)  node {$\times$};
\draw  (4.09434456222210, 2.20195121951220)  node {$\times$};
\draw  (7.04925484125584, 2.76873126873127)  node {$\times$};
\draw  (4.56434819146784, 2.56274570024570)  node {$\times$};
\draw  (2.77258872223978, 1.85700276789245)  node {$\times$};
\draw  (2.07944154167984, 0.825000000000000)  node {$\times$};
\draw  (2.89037175789616, 2.03529239521459)  node {$\times$};
\draw  (4.15888308335967, 1.99690865138085)  node {$\times$};
\draw  (2.30258509299405, 1.15588235294118)  node {$\times$};

    \end{tikzpicture}
  \end{small}
  \label{complex_graph}
\end{figure}

The dashed line has as equation $y = 5 ln(|G|)$. Therefore, we get
the following empirical overall complexity:
\begin{displaymath}
  \text{Computations} = O(ln(|G|) \times \text{Output size})
\end{displaymath}

\subsubsection{Tests around the unlabeled graph generation problem}

Although the generation engine is not optimized for the unlabeled
graph generation problem, we can apply our strategy on it.

Fix $n$, and consider the set $E$ of pairs of elements of $n$. The
symmetric group $\sg[n]$ acts on pairs by $\sigma \cdot (i,j) =
(\sigma(i), \sigma(j))$ for $\sigma \in \sg[n]$ and $(i,j)
\in E$. Let $G$ be the induced group of permutations of $E$. A
labeled graph can be identified with the integer vector with parts in
$0,1$. Then, two graphs are isomorphic if and only if the
corresponding vectors are in the same $G$-orbit. 

Now, one needs just to know which are these permutation groups
acting on pairs of integers. In the following example, we retrieve the
number of graphs on $n$ unlabeled nodes is, for small values of $n$ is
given by: $1$, $1$, $2$, $4$, $11$, $34$, $156$, $1044$, $12346$,
$274668$, $12005168$, ...

\sageex{
\sagepromt L = [TransitiveGroup(1,1), TransitiveGroup(3,2), \\
  TransitiveGroup(6,6), TransitiveGroup(10,12), TransitiveGroup(15,28), \\
  TransitiveGroup(21,38), TransitiveGroup(28,502)]

\sagepromt [IntegerVectorsModPermutationGroup(G,max\_part=1).cardinality() for G in L]

\sageret{[2, 4, 11, 34, 156, 1044, 12346]}
}


Notice that our generation engine generalizes the graph generation
problem in two directions. Removing the option {\tt max\_part}, one
enumerates multigraphs (graphs with multiple edges between
nodes). On the other hand, graphs correspond to special cases of
permutation groups. From an algebraic point of view, we saw graphs as
monomials whose exponents are $0$ or $1$, canonical for the action of
the symmetric group on pairs of nodes.

\section{Computing the invariants ring of a permutation group}

Let us explain how the generation engine from Section~\ref{enumeration_strategy}
is plugged into effective invariant theory
(see~\cite{kemper_derksen.cit.2002} and~\cite{king.2007.secondary}).

A well-known application to build an \emph{invariant polynomial} under
the action of a permutation group $G$ is the Reynolds operator
$R$. From any polynomial $P$ in $n$ variables $\x := x_1, x_2, \dots , x_n$,
the invariant is
\begin{displaymath}
  R(P) := \frac{1}{|G|} \sum_{\sigma \in G} \sigma\cdot P,
\end{displaymath}
where $\sigma\cdot P$ is the polynomial built from $P$ for which
$\sigma$ has permuted by position the tuple of variables $(x_1, x_2,
\dots , x_n)$. Formally, for any $\sigma \in G$
\begin{displaymath}
  (\sigma\cdot P)(x_1, x_2, ..., x_n) := P(x_{\sigma^{-1}(1)}, x_{\sigma^{-1}(2)}, \dots , x_{\sigma^{-1}(n)}).
\end{displaymath}

For large groups, the Reynolds operator is not very convenient to
build invariant polynomials. If $P$ is a monomial
$\x^\a := x_1^{a_1} x_2^{a_2} \cdots x_n^{a_n}$ where $\a = (a_1,
a_2,\dots , a_n)$, the minimal invariant one can build in
number of terms is the orbit sum $\displaystyle\orbsum[G](\x^{\a})$ of $\x$. 


Let $\KK$ a field, we denote by
$\KK[\mathbf{x}]^G$ the ring formed by all polynomials invariant under
the action of $G$.
\begin{displaymath}
  \KK[\x]^G := \{ P \in \KK[\x] | \forall \sigma \in G : \sigma \cdot P = P\}.
\end{displaymath}

For any subgroups $G$ of $\sg[n]$ and $\KK$ a field of characteristic
$0$, a result due to Hilbert and Noether state that the ring of
invariant $\KK[\x]^G$ is a free module of rank $\frac{n!}{|G|}$ over
the symmetric polynomials in the variable $\x$. Computing the
invariant ring $\KK[\x]^G$ consists essentially in building
algorithmically an explicit family (called \emph{secondary invariant
  polynomials}) of generators of this free module.

Searching the secondary invariant polynomials from orbit sum of
monomials whose vector of exponents is \emph{canonical} (instead of
all monomials) produces a gain of complexity of $|G|$ if we assume
that all orbits are of cardinality $|G|$. This assumption is obviously
false; however, in practice, it seems to hold in average and up to a
constant factor~\cite{Borie.2011.Thesis}).

In~\cite{Borie_Thiery.2011.Invariants}, the authors calculate the
secondary invariants of the $61^{st}$ transitive group over $14$
variables whose cardinality is $50 803 200$. Using the
\emph{canonicals} monomials, they managed to build a family of
$28$ irreducible secondary invariants deploying a set of $1716$
secondary invariants. This computation is unreachable by Gr\"obner
basis techniques.

\section{Computing primitive invariants for a permutation group}

\subsection{Introduction}

We now apply our generation strategy to this problem concerning
effective Galois theory.

\begin{problem}\label{pb_prim_invariant}
Let $n$ a positive integer and $G$ a permutation group, subgroup of
$\sg[n]$. Let $\KK$ be a field and $\x := x_1, \dots , x_n$ be $n$
formal variables. Find a polynomial $P \in \KK[x_1, \dots , x_n]$ such that
\begin{displaymath}
  \{ \sigma \in \sg[n] | \sigma \cdot P = P \} = G.
\end{displaymath}
A such polynomial is called a \emph{primitive invariant} for $G$.
\end{problem}

Problem~\ref{pb_prim_invariant} (exposed in~\cite{Girstmair.1987}
and~\cite{Abdeljaouad.TIATG}) consists in finding an invariant $P$
under the action of $G$ such that its stabilizer $Stab_{\sg[n]}(P)$ in
$\sg[n]$ is equal to $G$ itself. Solving this problem becomes difficult
when we want to construct a \emph{primitive invariant} of degree
minimal or a \emph{primitive invariant} with a minimal number of
terms.

\subsection{Primitive invariant of minimal degree}

\begin{small}
\begin{algorithm}[H]
\caption[]{Primitive invariant using stabilizer refinement}\label{stab_primitive_invariant}
  \raggedright
  
  \renewcommand{\labelitemi}{$\bullet$}
  Prerequisites : \\
  $\bullet$ $IntegerVectorsModPermgroup$: module to enumerate orbit
    representatives; \\
  $\bullet$ $stabilizer\_of\_orbit\_of(G, v)$: a function returning the
    permutation group which stabilizes the orbit of $v$ under the
    action of the permutation group $G$. 

  {\bf Arguments:} \\
  $\bullet$ $G$: A permutation group, subgroup of $\sg[n]$.

\begin{displaymath}
  \begin{array}{l}
    \textbf{def } minimal\_primitive\_invariant(G) : \\
    \ind cumulateStab \leftarrow SymmetricGroup(degree(G)) \\
    \ind chain \leftarrow [[(0,0, \dots, 0), cumulateStab, cumulateStab]] \\
    \ind \textbf{if } Cardinality(cumulateStab) == Cardinality(G): \\
    \ind \ind \textbf{return } chain \\
    \ind \textbf{for } v \in IntegerVectorsModPermgroup(G): \\
    \ind \ind AutV \leftarrow stabilizer\_of\_orbit\_of(G, v) \\
    \ind \ind Intersect \leftarrow cumulateStab \cap AutV \\
    \ind \ind \textbf{if } Cardinality(Intersect) < Cardinality(cumulateStab): \\
    \ind \ind \ind chain \leftarrow chain \cup [v, AutV, Intersect] \\
    \ind \ind \ind cumulateStab \leftarrow Intersect \\
    \ind \ind \ind \textbf{if } Cardinality(cumulateStab) == Cardinality(G): \\
    \ind \ind \ind \ind \textbf{return } chain \\
  \end{array}
\end{displaymath}
\end{algorithm}
\end{small}

\subsection{Benchmarks}

Algorithm~\ref{stab_primitive_invariant} terminates in less than an
hour for any subgroup of $\sg[10]$. Even, it can calculate some
primitive invariants for a lot of subgroups with degree between $10$
and $20$ while the literature only provides examples up to degree $7$
or $8$. Using the same computer, this benchmark just collects the average time
in seconds of execution of Algorithm~\ref{stab_primitive_invariant} by
executing systematically the algorithm on transitive groups of
degree $n$.

\begin{footnotesize}
\begin{table}[H]
  \centering
  \label{tab_call}
  \begin{tabular}{|c|c|c|c|c|c|c|c|c|c|}
    \hline
     Degree of Groups & 1 & 2 & 3 & 4 & 5 & 6 & 7 & 8 & 9 \\ \hline
     Computations time & 0.008 & 0.064 & 0.104 & 0.160 & 0.208 & 0.393 & 0.537 & 2.364 & 27.093 \\ \hline
  \end{tabular}
\end{table}
\end{footnotesize}

We would like to thanks Nicolas M. Thi\'ery, Simon A. King, Karl-Dieter
Crisman and Dmitri V. Pasechnik for useful comments about
implementation details, review of code and \cython optimizations.

This research was driven by computer exploration using the open-source
mathematical software \Sage~\cite{Sage}. In particular, we perused its
algebraic combinatorics features developed by the \sagecombinat
community~\cite{Sage-Combinat}, as well as its group theoretical
features provided by \gap~\cite{GAP}.

\bibliographystyle{alpha}
\bibliography{main}
\end{document}